\documentclass[a4paper, 11pt, twoside, final]{scrartcl}

\usepackage{PREAMBLE}
\usepackage{COMMAND}
\usepackage{HYPHENATION}

\begin{document}
\thispagestyle{empty}
\begin{center}
    \rule{\linewidth}{1pt}\\[0.4cm]
    {\sffamily \bfseries \large Non-negativity constraints in the one-dimensional\\
        discrete-time phase retrieval problem}\\[10pt] 
    {\sffamily\footnotesize Robert Beinert}\\[3pt]
    {\sffamily\footnotesize Institut für Numerische und Angewandte Mathematik}\\[-3pt]
    {\sffamily\footnotesize Georg-August-Universität Göttingen}\\
    \rule{\linewidth}{1pt}
\end{center}

\vspace*{10pt}

{\small
    \noindent
    {\sffamily\bfseries Abstract:} Phase retrieval problems occur in a width range of
    applications in physics and engineering such as crystallography, astronomy, and laser
    optics.  Common to all of them is the recovery of an unknown signal from the intensity
    of its Fourier transform.  Because of the well-known ambiguousness of these problems,
    the determination of the original signal is generally challenging.  Although there are
    many approaches in the literature to incorporate the assumption of non-negativity of
    the solution into numerical algorithms, theoretical considerations about the
    solvability with this constraint occur rarely.  In this paper, we consider the
    one-dimensional discrete-time setting and investigate whether the usually applied
    a~priori non-negativity can overcame the ambiguousness of the phase retrieval problem
    or not.  We show that the assumed non-negativity of the solution is usually not a
    sufficient a~priori condition to ensure uniqueness in one-dimensional phase retrieval.
    More precisely, using an appropriate characterization of the occurring ambiguities, we
    show that neither the uniqueness nor the ambiguousness are rare exceptions.

    \smallskip


    \smallskip
    
}

\section{Introduction}

In many fields of physics and engineering, one is faced with the recovery of an unknown
signal only from the intensity of its Fourier transform.  This phase retrieval problem
occurs in different applications as crystallography \cite{Mil90,Hau91}, astronomy
\cite{BS79,DF87} and laser optics \cite{SST04,SSD+06}.  In general, the recovery of an
analytic or numerical solution is challenging because of the well-known ambiguousness of
the problem.  Therefore, it is of essential importance to employ suitable a~priori
information about the original signal in order to find a meaningful solution or, in the
best case, the original signal itself.

In the rich literature on the phase retrieval problem, there are different approaches to
reduce the solution set or to ensure uniqueness.  For instance, the unknown signal $x$ can
be superposed with an appropriate reference signal $h$ such that one has access to the
additional Fourier intensity of $x+h$.  This idea was studied in \cite{KH90a,KH90,BP15}
for a known and in \cite{KH93,RDN13,BP15,Bei16} for an unknown reference signal.  More
particular reference signals have been considered in \cite{BFGR76,CESV13,Bei16}.  Instead
of interference measurements, it is also possible to use additional measurements in the
time domain.  For instance, one can employ additional magnitudes or phases to ensure a
unique recovery of the desired signal \cite{LT08,BP16}.

In the last years, the phase retrieval problem has been generalized from the classical
setting to the recovery of an finite-dimensional vector $x$ from appropriate frame
measurements $\absn{\iProd{x}{v_k}}$.  Here the question arises how the underlying frame
vectors have to be chosen, and how many frame vectors are needed to ensure the recovery of
$x$, see for instance \cite{BCE06, BBCE09,BCM14,BH15} and references therein.

In this paper, we consider the one-dimensional phase retrieval problem for discrete-time
signals, where we restrict ourselves to the recovery of an unknown signal with finite
support.  Here the occurring ambiguities can be explicitly specified by an appropriate
factorization of the autocorrelation signal, see \cite{BS79,BP15}.  Additionally, we
assume that the unknown signal is real-valued and non-negative.  This a~priori constraint
is usually applied if the unknown signal represents some intensity, see for instance
\cite{Fie78,BS79,DF87,SSD+06,LP14} and references therein.  Although there are many
efforts to incorporate the non-negativity into numerical algorithms, the solvability under
this constraint is studied rarely.  For this purpose, we consider the issue whether the
usually applied a~priori non-negativity can overcame the ambiguousness of the phase
retrieval problem or not.

The paper is organized as follows.  In section~\ref{sec:phase-retr-probl}, we introduce
the one-dimensional discrete-time phase retrieval problem and briefly recall the
characterization of the occurring ambiguities in \cite{BP15}.  Here we distinguish between
negligible, trivial ambiguities, like reflection and time shifts, and non-trivial
ambiguities.  Based on this characterization, we derive appropriate conditions whether a
solution is non-negative or not by exploiting that the Fourier transform of a
finite-supported signal is mainly an algebraic polynomial, see
section~\ref{sec:pos-algebr-polyn}.  Transferring our observation to the complete solution
set, we can explicitly construct phase retrieval problems that are uniquely solvable or
have a certain number of non-trivial non-negative solutions, see
section~\ref{sec:non-negat-sol}.  Finally, in section~\ref{sec:cont-prob-nonneg-sig}, we
present our main result that neither the ambiguousness nor the uniqueness are rare
exceptions, and that the non-negativity thus is not sufficient to ensure the unique
recovery of the desired signal.

\section{The phase retrieval problem}
\label{sec:phase-retr-probl}

In the following, we consider the \emph{one-dimensional discrete-time phase retrieval
    problem}.  This variant of the phase retrieval problem consists in the recovery of an
unknown discrete-time signal $x \coloneqq (x[n])_{n\in\Z}$ from its Fourier intensity
$\absn{\fourier x}$, where the discrete-time Fourier transform is given by
\begin{equation*}
    \fourier x \mleft( \omega \mright)
    \coloneqq \Fourier \mleft[ x \mright] \mleft( \omega \mright)
    \coloneqq \sum_{n\in\Z} x\mleft[n\mright] \, \e^{-\I \omega n}
    \qquad (\omega \in \R).
\end{equation*}
Further, we assume that the unknown signal $x$ has a finite support and that all
components $x[n]$ are non-negative.

Similarly to the recovery of a complex-valued signal, the phase retrieval problem for
non-negative signals always possesses some negligible ambiguities.  More precisely, we can
simply transfer \cite[Proposition~2.1]{BP15} to non-negative signals.

\begin{Proposition}\label{prop:trivial-amb}
    Let $x$ be a non-negative signal with finite support. Then
    \begin{enumerate}[\upshape(i)]
    \item\label{item:2} the time shifted signal $( x[n-n_0] )_{n\in\Z}$ for $n_0 \in \Z$
    \item\label{item:3} the reflected signal $(x[-n])_{n\in\Z}$
    \end{enumerate}
    have the same Fourier intensity $\absn{\fourier x}$.%
\end{Proposition}

Consequently, the applied assumption that the unknown signal is non-negative cannot ensure
uniqueness of the discrete-time phase retrieval problem.  However, since the shift and the
reflection in Proposition~\ref{prop:trivial-amb} are closely related to the original signal, we call
these negligible ambiguities \emph{trivial}.  Unfortunately, besides this trivial
ambiguities, our phase retrieval problem can have further non-trivial ambiguities as
exemplarily shown in \cite[Example~1, \etseqq]{BS79} and \cite[Figure~2]{Fie78}.  In order
to decide whether these examples are rare exceptions or the general case, we adapt the
characterization of the complete solution set in \cite{BP15} to our specific problem.

For this purpose, we recall that the \emph{autocorrelation signal} $a$ of a signal $x$ is
given by
\begin{equation*}
    a \mleft[ n \mright] \coloneqq
    \sum_{k\in\Z} \overline{x \mleft[ k \mright]} \, x \mleft[ k+n \mright]
    \qquad
    (n \in \Z),
\end{equation*}
and that the squared Fourier intensity can be written as
\begin{equation*}
    \absn{\fourier x \mleft( \omega \mright)}^2
    = \sum_{n\in \Z} \sum_{k \in \Z} x\mleft[n\mright] \, \overline{x\mleft[k\mright]} \, \e^{-\I\omega (n-k)}
    = \sum_{n\in \Z} \sum_{k \in \Z} x\mleft[k+n\mright] \, \overline{x\mleft[k\mright]} \, \e^{-\I\omega n}
    = \fourier a \mleft( \omega \mright) .
\end{equation*}
Since $x$ has a finite support, this property is transferred to the autocorrelation
signal.  Furthermore, the definition immediately implies that the components of $a$ have
to be symmetric, i.e., $a[-n] = a[n]$ for $n \in \Z$.  Thus, the \emph{autocorrelation
    function} $\fourier a$ is here always an even non-negative trigonometric polynomial of
degree $N-1$, where $N$ denotes the support length of the signal $x$.  Since a
trigonometric polynomial is completely determined by finitely many samples at appropriate
points, it is not necessary to know the Fourier intensity $\absn{\fourier x (\omega)}$ for
all $\omega \in \R$.  Indeed, the complete Fourier intensity of a real signal with support
length $N$ is already defined by $N$ samples in the interval $[0,\pi)$.

Following the lines in \cite{BP15}, we define the \emph{associated polynomial} $P$ to the
trigonometric polynomial $\fourier a$ by
\begin{equation*}
    P \mleft( z \mright)
    \coloneqq \sum_{n=0}^{\crampedclap{2N-2}} a \mleft[n-N+1 \mright] \, z^n
    \addmathskip
\end{equation*}
such that $\fourier a ( \omega ) = \e^{\I \omega (N-1)} \, P( \e^{-\I \omega} )$.  Since
the coefficients of $P$ are real and still satisfy $a[-n]=a[n]$, the zeros of $P$ have a
special structure.  More precisely, the real zeros occur in pairs
$(\gamma, \overline \gamma^{\,-1})$ and the complex zeros in quads
$(\gamma, \overline \gamma, \gamma^{-1}, \overline \gamma^{\,-1})$ or in the two pairs
$(\gamma,\overline \gamma^{\,-1})$ and $(\overline \gamma, \gamma^{-1})$. Thus, the
associated polynomial can always be written in the form
\begin{equation*}
    P \mleft( z \mright)
    = a \mleft[ N-1 \mright] \, \prod_{j=1}^{N-1} \left( z - \gamma_j^{\,} \right) \left(
        z - \overline \gamma_j^{\,-1} \right) \!.
    \addmathskip
\end{equation*}
Based on this observation, we can factorize the even non-negative polynomial $\fourier a$
by
\begin{align*}
  \fourier a \mleft( \omega \mright) 
  = \abs{ P \mleft( \e^{-\I \omega} \mright)}  
  &= \abs{a \mleft[ N-1 \mright]} \prod_{j=1}^{N-1} \abs{\e^{-\I\omega} -
    \gamma_j^{\,}} \abs{\e^{-\I\omega} - \overline \gamma_j^{\,-1}}
  \\
  &= \abs{a \mleft[ N-1 \mright]} \prod_{j=1}^{N-1} \abs{\e^{-\I\omega} -
    \gamma_j^{\,}} \absn{\gamma_j^{\,}}^{-1}  \, \absb{\overline \gamma_j - \e^{\I\omega}}
  \\
  &= \abs{ a \mleft[ N-1 \mright]} \prod_{j=1}^{N-1} \, \absn{\gamma_j}^{-1}
    \cdot \absB{ \prod_{j=1}^{N-1} \left( \e^{-\I \omega} - \gamma_j \right)}^2,
\end{align*}
which yields the following characterization of the solution set, see
\cite[Theorem~2.4]{BP15}.

\begin{Theorem}
    \label{the:repr-sol-time-dom}%
    Let $\fourier a$ be an even non-negative trigonometric polynomial of degree $N-1$.
    Then, each solution $x$ of the discrete-time phase retrieval problem
    $\absn{\fourier x}^{2} = \fourier a$ with finite support and non-negative components
    has a Fourier representation of the form
    \begin{equation}
        \label{eq:char-amb}
        \fourier x \mleft( \omega \mright)
        = \e^{-\I\omega n_0}
        \sqrt{\abs{a \mleft[ N-1 \mright]} \prod_{j=1}^{N-1} \abs{\beta_j}^{-1}} \cdot
        \prod_{j=1}^{N-1} \left( \e^{-\I \omega} - \beta_j \right) \!,
    \end{equation}
    where $n_0$ is an integer, and where for each $j$ the value
    $\beta_j$ is chosen from the zero pair $(\gamma_j^{\,}, \overline \gamma_j^{\,-1})$ of
    the associated polynomial to $\fourier a$.%
\end{Theorem}

Thus, each solution $x$ of the discrete-time phase retrieval problem
$\absn{\fourier x}^2 = \fourier a$ is uniquely given by the shift parameter $n_0$ and the
chosen values $\beta_j$.  Since $B \coloneqq \{\beta_1, \dots, \beta_{N-1}\}$ is a subset
of the zero set of the associated polynomial $P$, we call $B$ the \emph{corresponding zero
    set} of the solution $x$.  Besides the trivial shift ambiguity, which is directly
encoded in \eqref{eq:char-amb} by the factor $\e^{-\I \omega n_0}$,
Theorem~\ref{the:repr-sol-time-dom} covers the reflection ambiguity too.  More precisely,
one can show that the reflection $x[-\cdot]$ corresponds to the reflected zero set
\raisebox{0pt}[0pt][0pt]{$\{\overline \beta\kern0pt_1^{\,-1}, \dots, \overline
    \beta\kern0pt_{N-1}^{\,-1}\}$}
if $x$ corresponds to $\{\beta_1, \dots, \beta_{N-1}\}$.  Consequently, the discrete-time
phase retrieval problem to recover a non-negative signal $x$ with support length $N$ can
have at most $2^{N-2}$ non-trivially different solutions.

\section{Algebraic polynomials with non-negative coefficients}
\label{sec:pos-algebr-polyn}

To answer the question whether the phase retrieval problem in
Theorem~\ref{the:repr-sol-time-dom} can have more than one non-negative non-trivial
solution, we investigate conditions on the zero set
$B \coloneqq \{ \beta_1, \dots, \beta_{N-1} \}$ which ensure that a real signal with
finite support possesses only non-negative components.  We notice that the signal $x$ in
\eqref{eq:char-amb} is non-negative if and only if all coefficients of the monic
polynomial
\begin{equation*}
    Q \mleft( z \mright) \coloneqq \prod_{j=1}^{N-1} \left( z - \beta_j \right)
\end{equation*}
are non-negative.  Using Vieta's formulae and the \emph{elementary symmetric polynomials} $S_n$ defined by
\begin{equation*}
    S_n \mleft( \beta_1, \dots, \beta_{N-1} \mright)
    \coloneqq
    \smashoperator{\sum_{\quad1 \le k_1 < \dots < k_n \le N-1}} \beta_{k_1} \cdots
    \beta_{k_n}
    \qquad (n=1, \dots, N-1)
    \addmathskip
\end{equation*}
as well as $S_0 \coloneqq 1$ and $S_n \coloneqq 0$ for $n<0$ and $n \ge N$, we obtain the
representation
\begin{equation*}
    Q \mleft( z \mright)
    = \sum_{n=0}^{N-1} \left( -1 \right)^{n} \, S_{n} \mleft( \beta_1, \dots,
    \beta_{N-1} \mright) \, z^{N-1-n}.
\end{equation*}

The theorem of Descartes \cite[Satz~13.2]{Obr63} states that the number of positive zeros
of an algebraic polynomial with real coefficients is equal to the number of sign changes
in the coefficient sequence or less than it by an even number.  In our case, the
polynomial $Q$ has no sign changes, and thus all real zeros of the polynomial $Q$ have to
be negative.  In order to examine the dependency of the non-negativity of the coefficients
of $Q$ on the complex zero pairs, we generalize the observations in \cite{Bri85}.

\begin{Lemma}
    \label{lem:pos-cond-last-zero}
    Let $Q$ be a monic polynomial with real coefficients corresponding to the zero set
    $\{ \beta_1, \dots, \beta_{N-1} \}$. Assume that $(\beta_{N-2},\beta_{N-1})$ is a
    conjugated zero pair, and define
    $\sigma_n \coloneqq (-1)^n \, S_n (\beta_1, \dots, \beta_{N-3})$ for every $n \in
    \Z$. Then $Q$ has only non-negative coefficients if and only if $\beta_{N-1}$ fulfils
    \begin{equation}
        \label{eq:pos-cond-last-zero}
        \sigma_{n-2} \abs{ \beta_{N-1}}^2 - 2 \sigma_{n-1} \, \Re \beta_{N-1} + \sigma_n
        \ge 0
        \qquad(n=0, \dots, N-1).
    \end{equation}
\end{Lemma}

\begin{Proof}
    Since the zeros $\beta_{N-2}$ and $\beta_{N-1}$ form a conjugated pair, the monic
    polynomial $Q$ can be written as
    \vspace*{-10pt}
    \begin{equation*}
        Q \mleft( z \mright)
        = \left( z - \beta_{N-1}^{\,} \right) \left( z - \overline \beta_{N-1} \right)
        \prod_{j=1}^{N-3} \left( z - \beta_j \right) \!.
    \end{equation*}
    Observing that the product over the first $N-3$ linear factors is itself a monic
    polynomial, we can again apply Vieta's formulae and obtain
    \begin{align*}
        Q \mleft( z \mright)
        &= \left( z^2 - 2 \, \Re \beta_{N-1} \, z + \abs{\beta_{N-1}}^2 \right)
        \Bigl( \sum_{n=0}^{N-3} \sigma_n \, z^{N-3-n} \Bigr)
      \\[\fskip]
      &= \sum_{n=0}^{N-1} \left( \sigma_n - 2 \sigma_{n-1} \, \Re \beta_{N-1} +
            \sigma_{n-2} \abs{ \beta_{N-1}}^2 \right) z^{N-1-n},
    \end{align*}
     which completes the proof.  \qed%
\end{Proof}

\begin{Remark}
    \label{rem:pos-cond-last-zero}
    Each of the non-negativity constraints \eqref{eq:pos-cond-last-zero} describes a
    certain disc on the Riemann sphere.  This allows us to simplify the corresponding
    inequalities and to interpret them geometrically.  For example, if all zeros
    $\beta_1, \dots, \beta_{N-3}$ have a negative real part, one can show that the zero
    pair $(\beta_{N-1}, \overline \beta\kern0pt_{N-1}^{-1})$ has to lie in the closed half
    plane left of the imaginary axis through $\nicefrac{\sigma_1}{2}$ and, moreover, on or
    outside the circles with centre $ \nicefrac{\sigma_{n-1}}{\sigma_{n-2}}$ and radius
    \begin{equation*}
        \tfrac{\sqrt{\sigma_{n-1}^2- \sigma_n \sigma_{n-2}}}{\sigma_{n-2}}
        \qquad (n=2,\dots, N-2)
        \addmathskip
    \end{equation*}
    whenever the radius exists.  Indeed \eqref{eq:pos-cond-last-zero} implies $\Re
    \beta_{N-1} \le \nicefrac{\sigma_1}{2}$ for $n=1$ and
    \begin{equation*}
        \abs{\beta_{N-1} - \tfrac{\sigma_{n-1}}{\sigma_{n-2}}}^2  \ge
        \tfrac{\sigma_{n-1}^2 - \sigma_n \sigma_{n-2}}{\sigma_{n-2}^2}
    \end{equation*}
    for $n=2, \dots, N-2$.  This specific behavior is a complex version of the findings by
    Briggs in \cite[Section~7]{Bri85}. \qed
\end{Remark}

\section{Non-negative ambiguities of the phase retrieval problem}
\label{sec:non-negat-sol}

Based on our findings about the non-negativity of the coefficients of an algebraic
polynomial, we now investigate the non-negativity of the non-trivial solutions $x$ in
Theorem~\ref{the:repr-sol-time-dom}, which can be constructed by reflecting some of the
corresponding zeros $\beta_j$ at the unit circle.  First we show that, in the worst case,
the additional non-negativity constraint cannot reduce the set of non-trivial solutions at
all.

\begin{Proposition}
    \label{prop:nonneg-Hurw-signal}%
    Let $x$ be a real-valued discrete-time signal with finite support.  If the
    corresponding zero set $\{\beta_1, \dots, \beta_{N-1}\}$ is contained in the left half
    plane, \ie\ $\Re \beta_j < 0$ for all $j=1, \dots, N-1$, then all occurring
    real-valued non-trivial ambiguities of the corresponding phase retrieval problem are
    non-negative.%
\end{Proposition}

\begin{Proof}
    Using Theorem~\ref{the:repr-sol-time-dom}, we can generate all real-valued non-trivial
    ambiguities of the phase retrieval problem to recover $x$ by reflecting a subset of
    the real zeros \raisebox{0pt}[0pt][0pt]{$\beta_j$} and conjugate zero pairs
    $(\beta_j^{\,}, \overline \beta_j)$ at the unit circle.  Since all zeros $\beta_j$ and
    hence their reflections \raisebox{0pt}[0pt][0pt]{$\overline \beta\kern0pt_j^{\,-1}$}
    have a negative real part, the corresponding linear factors
    \begin{equation*}
        \e^{-\I\omega} - \beta_j
        \qquad \text{and} \qquad
        \left( \e^{-\I\omega} - \beta_j \right) \left( \e^{-\I\omega} - \overline \beta_j
        \right) = \e^{-2\I\omega} - 2 \, \Re[ \beta_j ] \, \e^{-\I\omega} + \absn{\beta_j}^2
    \end{equation*}
    of the real zeros $\beta_j$ and conjugate zero pairs
    $(\beta_j^{\,}, \overline \beta_j)$ in \eqref{eq:char-amb} have only non-negative
    coefficients. Thus, all possible non-trivial solutions have only non-negative
    components since a product of polynomials with non-negative coefficients has again
    non-negative coefficients.  \qed%
\end{Proof}

Besides this observation, we can exploit Theorem~\ref{lem:pos-cond-last-zero} to construct phase
retrieval problems with a specific number of non-negative non-trivial solutions.

\begin{Example}
    \label{ex:pos-cond-last-zero:sol-set}%
    We try to construct a phase retrieval problem with at least one non-negative solution
    $x$ by selecting the free conjugate zero pair $(\beta_4, \beta_5)$ of the
    corresponding zero set
    \begin{equation*}
        \Lambda \coloneqq \bigl\{ - \tfrac{3}{2}, -1+\I, -1-\I, \beta_4, \beta_5 \bigr\}
        \addmathskip
    \end{equation*}
    appropriately.  Since the reflection of the complete corresponding zero set leads to
    the reflection of the original signal, all further non-trivial solutions $y_1$, $y_2$,
    and $y_3$ according to Theorem~\ref{the:repr-sol-time-dom} are given by the zero sets
    \begin{align*}
        M_1 \coloneqq \bigl\{ - \tfrac{2}{3}, -1+\I, -1-\I, \beta_4, \beta_5  \bigr\},&
        &
        M_2 &\coloneqq \bigl\{ -\tfrac{3}{2}, -\tfrac{1}{2} \, ( 1 + \I ), - \tfrac{1}{2} \,
        (1 - \I ), \beta_4, \beta_5  \bigr\},
      \\[\fskip]
      \text{and}\qquad&&
        M_3 &\coloneqq \bigl\{ - \tfrac{2}{3}, -\tfrac{1}{2} \, ( 1 + \I ), - \tfrac{1}{2}
        \,  (1 - \I ), \beta_4, \beta_5  \bigr\},
    \end{align*}
    respectively.
  
    \begin{figure}[t]
        \centering
        \subfloat[{Restrictions on the last conjugate zero pair to ensure positivity of
            the signals $x$, $y_1$, $y_2$, and $y_3$
            \label{fig:pos-cond-last-zero:sol-set:restr}}]{%
            \makebox[5.7cm]{\includegraphics%
                {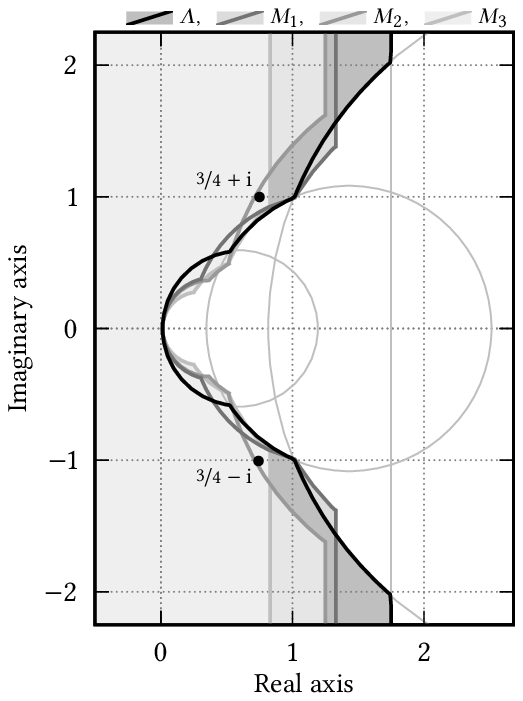}}}%
        \quad%
        \raisebox{80pt}{\parbox{5.8cm}{%
                \subfloat[{Signals $x$, $y_1$, $y_2$, and $y_3$ illustrated by
                    polygonal lines
                    \label{fig:pos-cond-last-zero:sol-set:sig}}]{%
                    \makebox[5.7cm]{\includegraphics%
                        {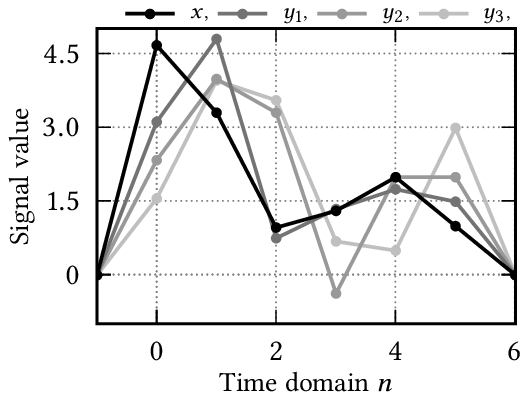}}}%
                \\%
                \subfloat[{Fourier intensities of the signals $x$, $y_1$, $y_2$, and
                    $y_3$ coincide
                    \label{fig:pos-cond-last-zero:sol-set:Fou}}]{%
                    \makebox[5.7cm]{\includegraphics%
                        {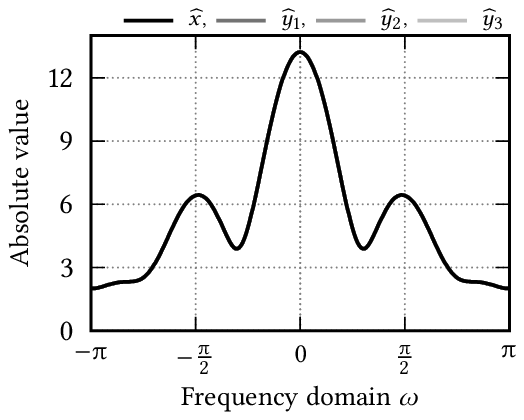}}}%
            }}%
        \caption{Restriction on the last zero pair in order to ensure non-negativity of
            different non-trivial ambiguities}
        \label{fig:pos-cond-last-zero:sol-set}
    \end{figure}

    The non-negativity constraints in Theorem~\ref{lem:pos-cond-last-zero} for these zero sets
    are visualized in Figure~\ref{fig:pos-cond-last-zero:sol-set}.  More detailed, the
    signal $x$ with zero set $\Lambda$ has only non-negative components if and only if the
    zero pair $(\beta_4, \beta_5)$ lies in the half plane left from the imaginary axis
    through $\nicefrac{7}{4}$ and outside the circles with centres $\nicefrac{7}{2}$,
    $\nicefrac{10}{7}$, and $\nicefrac{3}{5}$ in the complex plane and radii
    $\nicefrac{\sqrt{29}}{2}$, $\nicefrac{\sqrt{58}}{7}$, and $\nicefrac{3}{5}$
    respectively as discussed in Remark~\ref{rem:pos-cond-last-zero}.  The intersection of the
    half plane and the complements of the three discs is shown in
    Figure~\ref{fig:pos-cond-last-zero:sol-set}\subref{fig:pos-cond-last-zero:sol-set:restr}.
    The non-negativity constraints for the remaining sets $M_1$, $M_2$, and $M_3$ can be
    determined analogously.  Choosing $(\beta_4, \beta_5)$ in one or more intersections,
    we can thus ensure the non-negativity for certain non-trivial solutions and can
    directly influence the number of non-negative solutions.  For instance, if we choose
    \begin{equation*}
        \beta_4 \coloneqq \tfrac{3}{4} + \I
        \quad \text{and} \quad
        \beta_5 \coloneqq \tfrac{3}{4} - \I,
        \addmathskip
    \end{equation*}
    then the non-negativity constraints for $\Lambda$, $M_1$, and $M_3$ are fulfilled,
    which means that the phase retrieval problem to recover $x$ has two further
    non-negative non-trivial ambiguities.  The corresponding signals and Fourier
    intensities are shown in
    Figure~\ref{fig:pos-cond-last-zero:sol-set}\subref{fig:pos-cond-last-zero:sol-set:sig}
    and
    \ref{fig:pos-cond-last-zero:sol-set}\subref{fig:pos-cond-last-zero:sol-set:Fou}. \qed%
\end{Example}

\section{Uniqueness and ambiguousness under non-negativity constraints}
\label{sec:cont-prob-nonneg-sig}

Looking back at Example~\ref{ex:pos-cond-last-zero:sol-set}, it seems that the non-negativity
usually cannot reduce the number of arising non-trivial ambiguities.  However, the
situation dramatically depends on the fixed zeros $\beta_1, \dots, \beta_{N-3}$ in
Lemma~\ref{lem:pos-cond-last-zero}.  Although we cannot see the efficiency of the
non-negativity constraint directly, we can nevertheless use our findings to show that
neither uniqueness nor ambiguousness under the non-negativity constraint are rare
exceptions.  For this, we exploit that the non-trivial solutions continuously depend on
their corresponding zero sets, and vice versa.

\begin{Lemma}
    \label{lem:continuity-zeros}%
    Let $x$ be a discrete-time signal with support $\{0, \dots, N-1\}$ of length $N$ and
    corresponding zero set $\{\beta_1, \dots, \beta_{N-1} \}$.  For every sufficiently
    small number $\varepsilon > 0$, there exists a number $\delta > 0$ such that the
    corresponding zeros $\breve \beta_1, \dots, \breve \beta_{N-1}$ of every signal
    $\breve x$ with support $\{0,\dots, N-1\}$ of length $N$ and
    $\absn{ \breve x[n] - x[n]} \le \delta$ for $n$ from $0$ to $N-1$ can be ordered in a
    way that
    \begin{equation*}
        \abs{\breve \beta_j - \beta_j} \le \varepsilon
    \end{equation*}
    for $j$ from $1$ to $N-1$.%
\end{Lemma}

\begin{Proof}
    Based on the real-valued signal $x$, we consider the monic polynomial
    \begin{equation*}
        P\mleft( z \mright) = \frac{1}{x\mleft[N-1\mright]} \, \sum_{n=0}^{N-1} x \mleft[
        n \mright] \, z^n,
    \end{equation*}
    whose roots coincide with the zero set $\{\beta_1, \dots, \beta_{N-1} \}$.  In the
    following, we denote the coefficients of $P$ by
    $c_n \coloneqq \nicefrac{x[n]}{x[N-1]}$.  Using the continuity of roots theorem, see
    \cite[Theorem~3.1.1]{Ort72}, we find, for every sufficiently small number
    $\varepsilon > 0$, a number $\eta > 0$ such that the zeros $\breve \beta_j$ of all
    monic polynomials
    \begin{equation*}
        Q(z) \coloneqq z^{N-1} + \breve c_{N-2} \, z^{N-2} + \cdots + \breve c_0
    \end{equation*}
    with $\absn{ \breve c_n - c_n} \le \eta$ for $n$ from $0$ to $N-2$ can be ordered in 
    a way that
    \begin{equation*}
        \abs{\breve \beta_j - \beta_j} \le \varepsilon
        \submathskip
    \end{equation*}
    for $j$ from $1$ to $N-1$.

    If we identify $x$ with an $N$-dimensional vector, then the continuous mapping between
    the components $x[n]$ and the coefficients $c_n$ is given by
    \begin{equation*}
        \bigl( x\mleft[ 0 \mright], \dots, x \mleft[ N-1 \mright] \bigr)
        \mapsto
        \bigl( \tfrac{x\mleft[ 0 \mright]}{x\mleft[ N-1 \mright]}, \dots,
        \tfrac{x\mleft[ N-2 \mright]}{x\mleft[ N-1 \mright]} \bigr).
        \addmathskip
    \end{equation*}
    Hence, for every sufficiently small number $\eta>0$, there exists a number $\delta>0$
    such that the components of the image of every vector $\breve x$ in $\R^N$ with
    $\absn{\breve x[n] - x [n]} \le \delta$ for $n$ from $0$ to $N-1$ satisfy
    \begin{equation*}
        \abs{\tfrac{\breve x\mleft[ n \mright]}{\breve x \mleft[ N-1 \mright]}
            - \tfrac{ x\mleft[ n \mright]}{ x \mleft[ N-1 \mright]} } \le \eta
        \quad\text{or}\quad
        \abs{ \breve c_n - c_n} \le \eta
        \addmathskip
    \end{equation*}
    for $n$ from $0$ to $N-2$.  In order to avoid that $\breve x [N-1]$ becomes zero, we
    assume without loss of generality that $\delta < x[N-1]$.  Interpreting the vector
    $\breve x$ as discrete-time signal with support $\{0, \dots,N-1\}$ and combining both
    constructions yield the assertion.  \qed%
\end{Proof}

\begin{Remark}
    \label{rem:continuity-zeros}
    If we consider the discrete-time signals $\breve x$ in Lemma~\ref{lem:continuity-zeros} as
    $N$-di\-men\-sion\-al vectors, these signals form a closed ball with respect to the
    maximum norm.  Moreover, we can extend this ball to a cone since the multiplication of
    a signal with a positive real constant does not change the corresponding zero set.  By
    construction, the resulting cone cannot be contained in a set with zero Lebesgue
    measure. Consequently, this cone is an unbounded set with infinite measure.  \qed%
\end{Remark}

Applying Lemma~\ref{lem:continuity-zeros} to all possible non-trivial ambiguities in
Theorem~\ref{the:repr-sol-time-dom}, we can conclude that the occurring ambiguities continuously
depend on the original signal $x$.  In other words, for all signals $\breve x$ in a small
neighbourhood around $x$, the corresponding phase retrieval problem to recover $\breve x$
has the same solution behaviour as for $x$.  Moreover, this observation also holds for the
number of non-negative solutions if the corresponding zero set of $x$ fulfils some further
assumptions.  In the next two lemmata, we construct such signals.

\begin{Lemma}
    \label{lem:ass-neigh-amb-sol}%
    For every $N\in\N$, there exists a signal $x$ with support $\{0, \dots,N-1\}$ and
    positive components $x[n]$ for $n$ from $0$ to $N-1$ such that the phase retrieval
    problem to recover the signal $x$ has exactly $2^{N-2}$ non-trivial solutions
    satisfying the same assumptions.%
\end{Lemma}

\begin{Proof}
    We consider a signal $x$ with $N-1$ distinct real corresponding zeros fulfilling
    $\beta_j < -1$.  As a consequence, the signals in Theorem~\ref{the:repr-sol-time-dom} differ
    up to the trivial reflection ambiguity.  Choosing $n_0=0$, the phase retrieval problem
    to recover $x$ thus has $2^{N-2}$ non-trivially different solutions with support
    $\{0, \dots, N-1\}$.  The positivity of the components immediately follows from
    Proposition~\ref{prop:nonneg-Hurw-signal}.  \qed%
\end{Proof}

\begin{Lemma}
    \label{lem:ass-neigh-uni-sol}%
    For every $N > 3$, there exists a signal $x$ with support $\{0,\dots,N-1\}$, positive
    components $x[n]$ for $n$ from $0$ to $N-1$, and distinct zeros
    $\{\beta_1, \dots, \beta_{N-1}\}$ lying not on the unit circle such that the phase
    retrieval problem to recover the signal $x$ is uniquely solvable up to reflection.%
\end{Lemma}

\begin{Proof}
    Using the approach in Example~\ref{ex:pos-cond-last-zero:sol-set}, we choose distinct zeros
    $\beta_1, \dots, \beta_{N-3}$ with $\Re \beta_j < 1$ and extend this set by selecting
    an appropriate conjugate zero pair $(\beta_{N-2}, \beta_{N-1})$.  Since the fixed
    zeros $\beta_1, \dots, \beta_{N-3}$ lie in the left half plane, we can apply the
    slightly simpler constraints in Remark~\ref{rem:pos-cond-last-zero} to ensure the
    non-negativity of the corresponding signal $x$.  In this manner,
    $(\beta_{N-2}, \beta_{N-1})$ has to lie in the half plane left of the imaginary axis
    through $\nicefrac{\sigma_1}{2}$, which means that
    \begin{equation}
        \label{eq:half-plane-last-zero}
        \Re \beta_{N-1} \le
        - \tfrac{1}{2} \, \bigl( \Re \beta_1 + \cdots + \Re \beta_{N-3} \bigr) .
    \end{equation}
    By replacing a subset of zeros $\beta_j$ by their reflections at the unit circle, we
    obtain the non-negativity constraints for the remaining ambiguities in
    Theorem~\ref{the:repr-sol-time-dom} analogously.

    Since $\Re \beta_j < -1$ and thus $\Re \overline \beta\kern0pt_j^{\,-1} > -1$ for $j$
    from $1$ to $N-3$, the reflection of some zeros at the unit circle leads to a strictly
    smaller right-hand side of \eqref{eq:half-plane-last-zero}.  Consequently, we can
    choose $(\beta_{N-2}, \beta_{N-1})$ so that the zero set
    $\{\beta_1, \dots, \beta_{N-1} \}$ of the signal $x$ satisfies the non-negativity
    condition in \eqref{eq:half-plane-last-zero} strictly, and that the zero sets of the
    remaining non-trivial ambiguities violate this condition.  Figuratively, the zeros
    $\beta_{N-2}$ and $\beta_{N-1}$ have to lie in an appropriately small band in the
    complex plane, \cf\
    Figure~\ref{fig:pos-cond-last-zero:sol-set}\subref{fig:pos-cond-last-zero:sol-set:restr}.
    If we further ensure that the conjugate zero pair $(\beta_{N-2}, \beta_{N-1})$
    strictly lies outside the discs in Remark~\ref{rem:pos-cond-last-zero}, the constructed
    signal $x$ only possesses positive coefficients $x[n]$ for $n$ from $0$ to $N-1$ as
    desired. \qed%
\end{Proof}

We combine our findings in this section and finally show that neither the uniqueness nor
the ambiguousness under the non-negativity constraint is a rare exception.  Hence, the
assumed non-negativity of a discrete-time signal can be used to enforce the uniqueness of
the corresponding phase retrieval problem, but unfortunately not for every signal.

\begin{Theorem}
    \label{the:uni-amb-set-nonneg-sig}%
    The set of real-valued discrete-time signals with support $\{0,\dots,N-1\}$ of length
    $N>0$ that can be recovered uniquely up to reflection as well as the set of signals
    that cannot be recovered uniquely from their Fourier intensities employing the
    non-negativity constraint are both unbounded sets containing a cone of infinite
    Lebesgue measure.%
\end{Theorem}

\begin{Proof}
    In Lemma~\ref{lem:ass-neigh-amb-sol} and Lemma~\ref{lem:ass-neigh-uni-sol}, we have
    constructed signals $x$ so that the corresponding phase retrieval problem has either exactly
    $2^{N-2}$ non-trivial ambiguities or is uniquely solvable by choosing the
    corresponding zero set $\{\beta_1, \dots, \beta_{N-1}\}$ explicitly.  Since the
    non-negativity constraints in Lemma~\ref{lem:pos-cond-last-zero} and in
    Remark~\ref{rem:pos-cond-last-zero} continuously depend on the zeros $\beta_j$,
    there exists a small neighbourhood $U_B$ with respect to the maximum norm around the chosen
    set $B \coloneqq \{\beta_1, \dots, \beta_{N-1}\}$ such that the corresponding zero
    sets $\breve B \in U_B$ satisfy the same inequalities \eqref{eq:pos-cond-last-zero} as $B$.  Now,
    Lemma~\ref{lem:continuity-zeros} implies the existence of a small neighbourhood $U_x$
    around the constructed signal $x$ so that the corresponding phase retrieval problems
    have $2^{N-2}$ non-trivial solutions or are uniquely solvable respectively.  Extending
    the ball $U_x$ to a cone as discussed in Remark~\ref{rem:continuity-zeros} leads to the
    assertion.  \qed%
\end{Proof}

\section*{Acknowledgements}

I gratefully acknowledge the funding of this work by the DFG in the framework of the
SFB~755 \bq{Nanoscale photonic imaging} and of the GRK~2088 \bq{Discovering structure in
    complex data: statistics meets optimization and inverse problems.}

\bibliographystyle{alphadinUK}
{\footnotesize \bibliography{LITERATURE}}

\newcommand{\etalchar}[1]{$^{#1}$}
\begin{thebibliography}{BFGR76}


\providecommand{\url}[1]{\texttt{#1}}
\expandafter\ifx\csname urlstyle\endcsname\relax
  \providecommand{\doi}[1]{doi: #1}\else
  \providecommand{\doi}{doi: \begingroup \urlstyle{rm}\Url}\fi

\bibitem[BBCE09]{BBCE09}
\textsc{Balan}, Radu ; \textsc{Bodmann}, Bernhard~G. ; \textsc{Casazza},
  Peter~G.  ; \textsc{Edidin}, Dan:
\newblock Painless reconstruction from magnitudes of frame coefficients.
\newblock {In: }\emph{Journal of Fourier Analysis and Applications} 15 (2009),
  No. 4, pp. 488--501

\bibitem[BCE06]{BCE06}
\textsc{Balan}, Radu ; \textsc{Casazza}, Peter~G.  ; \textsc{Edidin}, Dan:
\newblock On signal reconstruction without phase.
\newblock {In: }\emph{Applied and Computational Harmonic Analysis} 20 (2006),
  No. 3, pp. 345--356

\bibitem[BCM14]{BCM14}
\textsc{Bandeira}, Afonso~S. ; \textsc{Chen}, Yutong  ; \textsc{Mixon},
  Dustin~G.:
\newblock Phase retrieval from power spectra of masked signals.
\newblock {In: }\emph{Information and Interference: A Journal of the IMA} 3
  (2014), June, No. 2, pp. 83--102

\bibitem[Bei16]{Bei16}
\textsc{Beinert}, Robert:
\newblock \emph{One-dimensional phase retrieval with additional interference
  measurements}.
\newblock April 2016. --
\newblock Preprint, arXiv:1604.04489v1

\bibitem[BFGR76]{BFGR76}
\textsc{Burge}, R.~E. ; \textsc{Fiddy}, M.~A. ; \textsc{Greenaway}, A.~H.  ;
  \textsc{Ross}, G.:
\newblock The phase problem.
\newblock {In: }\emph{Proceedings of the Royal Society of London. Series A.
  Mathematical Physical \& Engineering Sciences} 350 (1976), pp. 191--212

\bibitem[BH15]{BH15}
\textsc{Bodmann}, Bernhard~G. ; \textsc{Hammen}, Nathaniel:
\newblock Stable phase retrieval with low-redundancy frames.
\newblock {In: }\emph{Advances in Computational Mathematics} 41 (2015), April,
  No. 2, pp. 317--331

\bibitem[BP15]{BP15}
\textsc{Beinert}, Robert ; \textsc{Plonka}, Gerlind:
\newblock Ambiguities in one-dimensional discrete phase retrieval from Fourier
  magnitudes.
\newblock {In: }\emph{Journal of Fourier Analysis and Applications} 21 (2015),
  December, No. 6, pp. 1169--1198

\bibitem[BP16]{BP16}
\textsc{Beinert}, Robert ; \textsc{Plonka}, Gerlind:
\newblock \emph{Enforcing uniqueness in one-dimensional phase retrieval by
  additional signal information in time domain}.
\newblock March 2016. --
\newblock Preprint, arXiv:1604.04493v1

\bibitem[Bri85]{Bri85}
\textsc{Briggs}, William~E.:
\newblock Zeros and factors of polynomials with positive coefficients and
  protein-ligand binding.
\newblock {In: }\emph{Rocky Mountain Journal of Mathematics} 15 (1985), pp.
  75--89

\bibitem[BS79]{BS79}
\textsc{Bruck}, Yu.~M. ; \textsc{Sodin}, L.~G.:
\newblock On the ambiguity of the image reconstruction problem.
\newblock {In: }\emph{Optics communications} 30 (1979), September, No. 3, pp.
  304--308

\bibitem[CESV13]{CESV13}
\textsc{Cand{\`e}s}, Emmanuel~J. ; \textsc{Eldar}, Yonina~C. ;
  \textsc{Strohmer}, Thomas  ; \textsc{Voroninski}, Vladislav:
\newblock Phase retrieval via matrix completion.
\newblock {In: }\emph{SIAM Journal on Imaging Sciences} 6 (2013), No. 1, pp.
  199--225

\bibitem[DF87]{DF87}
\textsc{Dainty}, J.~C. ; \textsc{Fienup}, J.~R.:
\newblock Phase retrieval and image reconstruction for astronomy.
\newblock {In: }\textsc{Stark}, Henry (Ed.): \emph{Image Recovery {\upshape:}
  Theory and Application}.
\newblock Orlando (Florida) : Academic Press, 1987, Chapter ~7, pp. 231--275

\bibitem[Fie78]{Fie78}
\textsc{Fienup}, J.~R.:
\newblock Reconstruction of an object from the modulus of its Fourier
  transform.
\newblock {In: }\emph{Optics Letters} 3 (1978), July, No. 1, pp. 27--29

\bibitem[Hau91]{Hau91}
\textsc{Hauptman}, Herbert~A.:
\newblock The phase problem of x-ray crystallography.
\newblock {In: }\emph{Reports on Progress in Physics} 54 (1991), November, No.
  11, pp. 1427--1454

\bibitem[KH90a]{KH90a}
\textsc{Kim}, Wooshik ; \textsc{Hayes}, Monson~H.:
\newblock Iterative phase retrieval using two Fourier transform intensities.
\newblock {In: }\emph{Proceedings {\upshape:} ICASSP~90 {\upshape:} 1990
  International Conference on Acoustics, Speech and Signal Processing
  {\upshape:} April 3--6, 1990} Vol.~3 IEEE Signal Processing Society, 1990,
  pp. 1563--1566

\bibitem[KH90b]{KH90}
\textsc{Kim}, Wooshik ; \textsc{Hayes}, Monson~H.:
\newblock Phase retrieval using two Fourier-transform intensities.
\newblock {In: }\emph{Journal of the Optical Society of America A} 7 (1990),
  March, No. 3, pp. 441--449

\bibitem[KH93]{KH93}
\textsc{Kim}, Wooshik ; \textsc{Hayes}, Monson~H.:
\newblock Phase retrieval using a window function.
\newblock {In: }\emph{IEEE Transactions on Signal Processing} 41 (1993), March,
  No. 3, pp. 1409--1412

\bibitem[LP14]{LP14}
\textsc{Loock}, Stefan ; \textsc{Plonka}, Gerlind:
\newblock Phase retrieval for Fresnel measurements using a shearlet sparsity
  constraint.
\newblock {In: }\emph{Inverse Problems} 30 (2014), No. 5, pp. 055005(17)

\bibitem[LT08]{LT08}
\textsc{Langemann}, Dirk ; \textsc{Tasche}, Manfred:
\newblock Phase reconstruction by a multilevel iteratively regularized
  {G}auss-{N}ewton method.
\newblock {In: }\emph{Inverse Problems} 24 (2008), No. 3, pp. 035006(26)

\bibitem[Mil90]{Mil90}
\textsc{Millane}, R.~P.:
\newblock Phase retrieval in crystallography and optics.
\newblock {In: }\emph{Journal of the Optical Society of America A} 7 (1990),
  March, No. 3, pp. 394--411

\bibitem[Obr63]{Obr63}
\textsc{Obreschkoff}, Nikola:
\newblock \emph{Verteilung und Berechnung der Nullstellen reeller Polynome.}
\newblock Berlin : VEB Deutscher Verlag der Wissenschaften, 1963
  (Hochschulb{\"u}cher f{\"u}r Ma\-the\-ma\-tik~55)

\bibitem[Ort72]{Ort72}
\textsc{Ortega}, James~M.:
\newblock \emph{Numerical Analysis {\upshape:} A Second Course}.
\newblock New York : Academic Press, 1972 (Computer Science and Applied
  Mathematics)

\bibitem[RDN13]{RDN13}
\textsc{Raz}, Oren ; \textsc{Dudovich}, Nirit  ; \textsc{Nadler}, Boaz:
\newblock Vectorial phase retrieval of 1-D signals.
\newblock {In: }\emph{IEEE Transactions on Signal Processing} 61 (2013), April,
  No. 7, pp. 1632--1643

\bibitem[SSD{\etalchar{+}}06]{SSD+06}
\textsc{Seifert}, Birger ; \textsc{Stolz}, Heinrich ; \textsc{Donatelli}, Marco
  ; \textsc{Langemann}, Dirk  ; \textsc{Tasche}, Manfred:
\newblock Multilevel {G}auss-{N}ewton methods for phase retrieval problems.
\newblock {In: }\emph{Journal of Physics. A. Mathematical and General} 39
  (2006), No. 16, pp. 4191--4206

\bibitem[SST04]{SST04}
\textsc{Seifert}, Birger ; \textsc{Stolz}, Heinrich  ; \textsc{Tasche},
  Manfred:
\newblock Nontrivial ambiguities for blind frequency-resolved optical gating
  and the problem of uniqueness.
\newblock {In: }\emph{Journal of the Optical Society of America B} 21 (2004),
  May, No. 5, pp. 1089--1097

\end{thebibliography}

\end{document}